\documentclass[12 pt]{amsart}
\usepackage{amscd,amssymb}
\usepackage[arrow,matrix]{xy}
\usepackage{graphicx}
\oddsidemargin=0.3in \evensidemargin=0.3in \theoremstyle{plain}
\newtheorem{theorem}[subsection]{Theorem}
\newtheorem{lemma}[subsection]{Lemma}

\newtheorem{cor}[subsection]{Corollary}

\theoremstyle{definition}
\newtheorem{rk}[subsection]{Remark}

\numberwithin{equation}{section} 

\newcommand{\RR}{{\mathcal R}}

\newcommand{\A}{{\mathcal A}}
\newcommand{\B}{{\mathcal B}}

\newcommand{\LL}{{\mathcal L}}

\newcommand{\V}{{\mathcal V}}

\newcommand{\Z}{\mathbb{Z}}
\newcommand{\Q}{\mathbb{Q}}

\newcommand{\C}{\mathbb{C}}

\newcommand{\T}{\mathbb{T}}

\DeclareMathOperator{\Hom}{Hom} \DeclareMathOperator{\rank}{rank}
 
\begin{document}
\title [On the intersection of rational transversal subtori]
{On the intersection of rational transversal subtori }
\author{SHAHEEN NAZIR}
 \address{School of Mathematical Sciences,
         Government College University,
         68-B New Muslim Town Lahore,
         PAKISTAN.}
\email {shaheen.nazeer@gmail.com}
\subjclass[2000]{Primary 14C21,
14F99, 32S22 ; Secondary 14E05, 14H50.}

\keywords{local system, constructible sheaf, twisted cohomology,
characteristic variety, pencil of plane curves}

\begin{abstract}
 We show that under a suitable transversality condition, the
intersection of two rational subtori in an algebraic torus
$(\C^*)^n$ is a finite group which can be determined using the
torsion part of some associated lattice. Applications are given to
the study of characteristic varieties of smooth complex algebraic
varieties. As an example we discuss A. Suciu's line arrangement, the
so-called deleted $B_3$-arrangement.

\end{abstract}

\maketitle

\section{Introduction} \label{sec:intro}

Let $L$ be a free $\Z$-module of finite rank $n$, and let $A \subset
L$ and $B \subset L$ be two primitive sublattices, i.e. $A$ and $B$
are subgroups such that $Tors (L/A)=Tors (B/A)=0.$ Consider the
associated $\C$-vector spaces
$$V=L \otimes _{\Z}\C, V_A=A \otimes _{\Z}\C \text{  and  } V_B=B \otimes _{\Z}
\C.$$ Let $\exp_L: V \to T=L \otimes _{\Z}\C^*$ be the associated
exponential map given by
$$\exp_L=1_L \otimes _{\Z} \exp$$
where $1_L:L \to L$ is the identity and $\exp:\C \to \C^*$ is
defined by $t \mapsto \exp(2\pi it)$. Then $\exp_L$ is a surjective
group homomorphism with kernel $L=L\otimes _{\Z}\Z \subset V$. If a
$\Z$-basis of $L$ is chosen, then one has obvious identifications
$L=\Z^n$, $V=\C^n$, $T=(\C^*)^n$ and $\exp_L: \C^n \to (\C^*)^n$ is
given by
$$(t_1,...,t_n) \mapsto (\exp(2\pi i t_1),...,\exp(2\pi i t_n)).$$
The main result of this note is the following.

\begin{theorem} \label{theorem1}
With the above notation, if in addition $V_A \cap V_B=0$, then there
is a group isomorphism
$$Tors(L/(A+B)) \to \exp_L(V_A) \cap \exp_L(V_B).$$
\end{theorem}

In fact any algebraic subtorus $S \subset T$, (i.e. $S$ is a closed
algebraic subset and a subgroup in $T$), comes from a primitive lattice
$A(S) \subset L$, see Lemma 2.1 in Section II in Arapura's paper
\cite{A}.
Hence Theorem \ref{theorem1} applies to any pair of such algebraic subtori.

This Theorem is proved in the second section. In the third section
we show how to use Theorem \ref{theorem1} to describe the
intersections of the irreducible components of the characteristic
varieties of smooth complex algebraic varieties. A specific example
coming from the hyperplane arrangement theory concludes the paper.

\section{The proof} \label{sec:proof}

Let $n=\rank L$, $a=\rank A$ and $b=\rank B$. Consider the quotient
group $L'=L/A$, which is again a lattice, of rank $n-a$. The
composition $B \to L \to L'$ of the inclusion $B \to L$ and the
projection $L \to L'$ gives rise to an injective morphism $\iota: B
\to L'$ identifying $B$ to the sublattice $B'=\iota (B) \subset L'$.

Then there is a basis $e'_1,...,e'_{n-a}$ of the lattice $L'$ such
that $B'$ is the subgroup spanned by $d_1e'_1,...,d_be'_b$ for some
positive integers $d_j$. Moreover there is an integer $m$ with $1
\le m \le b+1$ such that
\begin{equation} \label{eq1}
1=d_1=...=d_{m-1}<d_m \le ... \le d_b \text{ and }
d_m|d_{m+1}|...|d_b.
\end{equation}
It follows that
\begin{equation} \label{eq2}
Tors (L/(A+B))=Tors(\frac{L/A}{(A+B)/A})=Tors(L'/B')
\end{equation}
and hence
\begin{equation} \label{eq21}
Tors (L/(A+B))=\Z/ d_m\Z \oplus \Z/d_{m+1}\Z \oplus ... \oplus
\Z/d_b\Z.
\end{equation}

Let $e_1,...,e_{n-a}$ be any lifts of the vectors $e'_j$'s to $L$
and let $f_1,...,f_a$ be a $\Z$-basis of $A$. Then
$\B=\{e_1,...,e_{n-a},f_1,...,f_a\}$ is a $\Z$-basis of $L$.

For $j=1,...,b$, let $g_j \in B$ be vectors such that their classes
$g'_j$ in $L'$ satisfy $g'_j=d_j e'_j$. It follows that
$g_j=d_je_j+a_j$ for some vectors $a_j \in A$. Write now
\begin{equation} \label{eq3}
a_j=\sum_{i=1,a}\alpha _{ji}f_i
\end{equation}
for some $\alpha _{ji} \in \Z$. By replacing $e_j$ by $e_j+r_j$ for
suitable vectors $r_j \in A$, we may and will assume in the sequel
that
\begin{equation} \label{eq4}
0\le \alpha _{ji} <d_j
\end{equation}
for all $i=1,...,a$ and $j=1,...,b$. In particular $a_j=0$ for
$j=1,...,m-1$.

\begin{lemma} \label{lemma1}
The vectors $g_1,...,g_b$ form a $\Z$-basis of the lattice $B$.

 \end{lemma}

\proof Note that the vectors $g_1,...,g_b$ are all contained in $B$
and, on the other hand, their images under $\iota$ span the lattice
$B'$.

\endproof

Assume now that $\exp_L(v_A)=\exp_L(v_B)$ for some vectors
$$v_A=p_1f_1+...+p_af_a \in V_A$$
and
$$v_B=q_1g_1+...+q_bg_b \in V_B$$
where $p_i,q_j \in \C$. It follows that $v_A-v_B \in \ker \exp_L=L$.
More precisely, we get
$$q_jd_j \in \Z \text{ for } j=1,...,b$$
and
$$z_i:=p_i-\sum _{j=1,b}q_j\alpha _{ji} \in \Z \text{ for } i=1,...,a.$$
It follows that $q_j=k_j/d_j$ and we may and will assume that $0\le
k_j<d_j$, since the value of $\exp_L(v_B)$ is not changed when the
coefficients $q_j$ are modified by integers. Note that with this
choice one has $k_j=0$ for $j=1,...,m- 1$. In this way we get a
surjective group homomorphism
\begin{equation} \label{eq5}
\theta: \Z/ d_m\Z \oplus \Z/d_{m+1}\Z \oplus ... \oplus \Z/d_b\Z \to
\exp_L (V_A) \cap \exp_L(V_B)
\end{equation}
given by
$$\hat k =(\hat k_m,...,\hat k_b) \mapsto \exp_L(\frac{k_m}{d_m}g_m+...+\frac
{k_b}{d_b}g_b).$$ This morphism $\theta$ is indeed correctly defined
since for any choice of the $q_j$'s as above we may use the defining
equation of $z_i$ above, set $z_i=0$ and determine the values for
$p_i$'s, i.e. find a vector $v_A$ such that
$\exp_L(v_A)=\exp_L(v_B)$.

To show that $\theta$ is injective, we have to show that $\ker
\theta=0$.

Since $B$ is primitive so on the same lines we can take the set\\
$\{g_1,\cdots,g_b,h_1,\cdots,h_{n-b}\}$ as a $\Z-$basis of $L$,
 where $h_1',\cdots,h_{n-b}'$ is a $\Z-$basis for the lattice $L/B$.
 Let $\hat{k}\in \ker \theta.$ Then
 $\theta(\hat{k})=\exp_L(\frac{k_m}{d_m}g_m+\cdots+\frac{k_b}{d_b}g_b)=1$,
 which implies that $\frac{k_i}{d_i}\in \Z,$ for all $m\leq i \leq
 b.$ Therefore, $\hat{k}=(\hat k_m,...,\hat k_b)=(\hat 0,...,\hat 0)$ i.e., $\ker\theta=0.$

\section{On the intersection of irreducible components of characteristic
varieties} \label{sec:char var}

\subsection{Local systems, characteristic and resonance varieties} \label
{sec:3.2}

Let $M$ be a quasi-projective smooth complex algebraic variety. The
rank one local systems on $M$ are parameterized by the algebraic
group
\begin{equation} \label{eq9}
\T(M)= \Hom( H_1(M),\C^*)
\end{equation}
The connected component $\T^0(M)$ of the unit element $1 \in \T(M)$
is an algebraic torus, i.e. it is isomorphic to $(\C^*)^n$, where $n
\in\mathbb{N}$ is the first batti number of $M$, i.e., $n=b_1(M)$.
It is clear that $\T^0 (M)=\T(M)$ if and only if the integral
homology group $H_1(M)$ is torsion free. For $\rho \in \T(M)$, we
denote by $\LL_{\rho}$ the corresponding local system on $M$.

The computation of the twisted cohomology groups $H^j(M, \LL_{
\rho})$ is one of the major problems in many areas of topology. To
study these cohomology groups, one idea is to study the {\it
characteristic varieties} defined by
\begin{equation} \label{eq10}
\V_m^j(M)=\{ \rho \in \T(M) ~|~  \dim H^j(M, \LL_{ \rho}) \ge m \}.
\end{equation}
To simplify the notation, we set $\V_m(M)=\V_m^1(M).  $
It is known that the following holds, see Beauville \cite{Beau} and
Simpson \cite{Sim} in the proper case and Arapura \cite{A} in the
quasi-projective case.

\begin{theorem} \label{theorem2}
The positive dimensional irreducible components of $\V_m(M)$ are
translated subtori in $\T(M)$ by elements of finite order. More
precisely, for each positive dimensional irreducible component $W$
of $\V_m(M)$  we can write $W= \rho \cdot f^*(\T(S))$, where $f:M
\to S$ is a surjective regular mapping to a curve $S$, having a
connected general fiber and $\rho \in \T(M)$ is a finite order
character.
\end{theorem}

If $1 \in W$, then we can take $\rho =1$ in the above equality. Let
$T_1W$ denote the tangent space to $W$ at the identity $1$ in such a
case. Theorem 2, (b) in  \cite{DPS} gives the following.

\begin{theorem} \label{theorem3}
Let $M$ be a  quasi-projective smooth complex algebraic variety. Let
$W_1$ and $W_2$ be two distinct irreducible components of the
characteristic variety $\V_1 (M)$ such that $1 \in W_1 \cap W_2$.
Then $T_1W_1\cap T_1W_2=0$.
\end{theorem}

Note  that any such  tangent space $T_1W \subset H^1(M,\C)$ is {\it
rationally defined}, i.e. there is a primitive lattice $L \subset
H^1(M,\Z)$ such that $T_1W =L \otimes _{\Z}\C$ under the
identification $ H^1(M,\C)= H^1 (M,\Z) \otimes _{\Z}\C$. Indeed, one
can take $L$ to be the primitive sublattice $f^*(H^1(S,\Z))$, in
view of the functoriality of the exponential mapping $\exp:
T_1\T(M)=H^1(M,\C) \to \T(M)$ and of the following.

\begin{lemma} \label{lemma2}
 Let $f:M \to S$ be a surjective regular mapping to a curve $S$, having a
connected general fiber.
 Then $f^*(H^1(S,\Z))$ is a primitive sublattice in $H^1(M,\Z)$.

 \end{lemma}

\proof In these conditions, it is well known that the morphism
$f_*:H_1(M,\Z) \to H_1 (S,\Z)$ is surjective. Let $L_0$ be a
primitive sublattice in $H_1(M,\Z)$ such that $H_1(M,\Z)=\ker f_*
\oplus L_0$. Then $f^*(H^1(S,\Z))$ can be identified to the dual
$L_0^{\vee}=\Hom(L_0,\Z) =\{u \in  \Hom(H_1(M,\Z),\Z) ~~;  ~~ u|\ker
f_*=0\} $. Moreover $H^1(M,\Z)=H_1(M,\Z)^{\vee}= (\ker f_*)^{\vee}
\oplus L_0^{\vee} $, which completes the proof of the claim.

\endproof

Applying Theorem \ref{theorem1} to this setting, we get the
following.

\begin{cor} \label{cor1}
Let $W_1$ and $W_2$ be two distinct irreducible components of the
characteristic variety $\V_1(M)$ such that $1 \in W_1 \cap W_2$. Let
$L_1$ and $L_2$ be the primitive sublattices in $ H^1(M,\Z)$
associated to $W_1$ and $W_2 $ respectively by the above
construction. Then there is a group isomorphism
$$Tors(H^1(M,\Z)/(L_1+L_2)) =W_1 \cap W_2.$$
\end{cor}

\begin{rk} \label{rk1}
Let $W_1$ and $W_2$ be two distinct irreducible components of the
characteristic variety $\V_1(M)$, at least one of them, say $W_1$
translated, i.e. $1 \notin W_1$, and meeting at a point $\rho$. Then
we may write $W_1=\rho \cdot W_1'$ and $W_2=\rho \cdot W_2'$, where
$W_j'$ are subtori in $\T(M)$.

Assume that $\dim W_j >1$ for $j=1,2$ (the claim is obvious when one
of the two components is one dimensional) and that $M$ is a
hyperplane arrangement complement. Then  $T_1W_1' \cap T_1W_2'=0$,
since $W'_1$ and $W'_2$ are again two distinct irreducible
components of the characteristic variety $\V_1(M)$, see  \cite{D1}.
Moreover, one can see exactly as above that the tangent spaces
$T_1W_j'$ are rationally defined and the above Corollary yields a
set bijection
$$Tors(H^1(M,\Z)/(L_1'+L_2')) =W_1 \cap W_2$$
where $L_j'$ is the primitive sublattice associated to $W_j'$ by the
above construction.

Note that any character in such an intersection $W_1 \cap W_2$ has
finite order. Indeed, let $W_1=\rho_1.W_1'$ and $W_2=\rho_2.W_2'$,
where $\rho_1 $ and $\rho_2$ are finite order characters, i.e.
$\rho_1^{m_1}=1$ and $\rho_2^{m_2}=1$. Then $\rho\in W_1\cap W_2$
implies that $\rho=\rho_1 w_1=\rho_2 w_2,$ where $w_j\in W_j'$. Let
$m=lcm(m_1,m_2)$. Then $\rho^m=\rho_1^m w_1^m$ and $\rho^m=\rho_2^m
w_2^m$ which implies that $\rho^m=w_1^m=w_1^m\ \Rightarrow$
$\rho^m\in W'_1\cap W'_2$ so by Corollary  \ref{cor1}, $\rho^m$ is
of finite order. Thus, $\rho$ is a finite order character.

A completly different proof of the finiteness of the intersection
 $W_1 \cap W_2$ of two distinct irreducible components of the first
 characteristic variety was given in \cite{DPS3}.

\end{rk}

Let $H^*(M,\C)$ be the cohomology algebra of the variety $M$ with
$\C$- coefficients. Right multiplication (cup-product) by an element
$z \in H^1(M,\C) $ yields a cochain complex $(H^*(M,\C), \mu_z)$.
The {\it resonance varieties} of $M$ are the jumping loci for the
 cohomology of this complex, namely
\begin{equation} \label{rv1}
\RR_m^j(M)=\{z \in H^1(M,\C) ~|~ \dim H^j(H^*(M,\C), \mu_z) \ge m \}.
\end{equation}
To simplify the notation, we set $\RR_m(M)=\RR_m^1(M)$

One of the main results relating the characteristic and resonance
varieties is the following.

\begin{theorem} \label{thm2}
Let $M$ be a hypersurface arrangement complement. The exponential
map $\exp: H^1(M,\C) \to \T^0(M)$ induces for any $m,j \ge 1$ an
isomorphism of analytic germs
$$(\RR_m^j(M),0) \simeq (\V_m^j(M),1).$$

\end{theorem}

This equality of germs implies that the resonance variety $\RR_m^j(M)$
is exactly the tangent cone at 1 of the characteristic variety
$\V_m^j(M)$, a fact established by Cohen and Suciu \cite{CS1} and which
can also be derived from  \cite{ESV}. See also Theorem 3.7 in
\cite{CO}. It was
claimed by A. Libgober that this property holds for any smooth
quasi-projective variety, see \cite{Li}, but now there are
counter-examples to this claim, see \cite{DPS}.

\begin{rk} \label{rk13}
 According to Theorem 1.1 in Section V in Arapura \cite{A}, under the
 assumption that $H^1(M,\Q)$
has a pure Hodge structure, the positive dimensional irreducible
 components
 of all characteristic varieties $\V_m^j(M)$ are (translated) subtori.
 Our Theorem \ref{theorem1} applies to this more general setting as
 well. The major difference with the case of first characteristic
 varieties
 is that distinct components do not necessarily meet only at the origin.
 Here is an example for which we are grateful to Professor A. Suciu.

 Consider the central hyperplane arrangement in $\C^4$ defined by the equation
$$xyzw(x+y+z)(y-z+w)=0.$$
 Then the corresponding resonance variety  $\RR_1^2(M)$
 has two 3-dimensional components, say $E_1$ and $E_2$, defined
 respectively by the ideals
 $$I_1=(x_1+x_2+x_3+x_6,x_4,x_5)$$
 and
 $$I_2=(x_1, x_2+x_3+x_4+x_5,x_6).$$
 These two components intersect in the line
 $D=(x_1,x_2+x_3,x_4,x_5,x_6)$.
 This implies that the irreducible components $W_1=\exp(E_1)$
 and  $W_2=\exp(E_2)$ of the characteristic variety  $\V_1^2(M)$
 intersect along  the 1-dimensional subtorus $\exp(D)$.

\end{rk}

The fact that $M$ has a simply-connected compactification implies
the following.

\begin{cor} \label{cor0.1}
The irreducible components of $\RR_1(M)$ are precisely the maximal
linear subspaces $E \subset H^1(M,\C)$, isotropic with respect to
the cup product on $M$
$$\cup: H^1(M,\C) \times H^1(M,\C) \to H^2(M,\C)$$
 and such that $\dim E \ge 2$.
\end{cor}

\proof

Let $E$ be a component of $\RR_1(M)$. By the above Theorem there is
a component $W$ in $\V_1(M)$ such that $1 \in W$ and $T_1W=E$. By
Arapura's results in \cite{A} we can write $W=f_E^*(\T(S))$, where
$f_E:M \to S$ is a regular mapping to a curve $S$. Since in our case
$S$ is rational, $T_1W=f_E^*(H^1(S,\C))$ is isotropic
 with respect to the cup product, since the cup product on $H^1(S,\C)$ is trivial.
Maximality of $E$ comes from the fact that $E$ is a component of
$\RR_1(M)$. The restriction $\dim E \ge 2$ comes from  \cite{A}. A
mapping $f_E$ as above is said to be associated to the subspace $E$.

\endproof

\section{An example: the deleted $B_3$-arrangement} \label{sec:3.3}

Let $\A$ be the deleted $B_3$-arrangement which is obtained from the
$B_3$ reflection arrangement by deleting the plane $x+y-z=0.$ A
defining polynomial for $\A$ is
$Q=xyz(x-y)(x-z)(y-z)(x-y-z)(x-y+z)$. The decone $d\A$ is obtained
by setting $z=1$. Let $L_1:\ell_1=x=0,\ L_2:\ell_2=y=0,\
L_3:\ell_3=x-y=0,\ L_4:\ell_4=x-1=0,\ L_5: \ell_5=y-1=0,\ L_6:
\ell_6=x-y-1=0,\ L_7: \ell_7=x-y+1=0$ be the lines of the associated
affine arrangement in $\mathbb{A}^2.$ Let $L_8: \ell_8=z=0$
 be the line at infinity and let $M$ be the complement of $d\A$ in $\mathbb{A}^2$.
The resonant variety $R_1(d\A)$ has 12 irreducible components of
dimension 2 and 3. These components $E$ and their associated maps
$f_E: M \rightarrow S$ are given below, see \cite{S1}, \cite{S2}.
Denote by $e_1,...,e_7$ the $\Z$-basis of $H^1(M,\Z)$ given by
$e_j=\frac {1}{2\pi i}\frac{d\ell_j}{\ell_j}$, see
 \cite{OT}. Then each of the components $E$ is the $\C$-vector space spanned by a primitive lattice denoted by $E ^0$, i.e. $E=E^0 \otimes _{\Z}\C$, so it is enough in each case to indicate a $\Z$-basis of $E^0$.

\medskip

{1. \bf The local components} There are 7 local components,
corresponding to 6 triple points and 1 quadruple point.\\
For the triple $L_1\cap L_2\cap L_3$, we have\\
$E_1^0=<e_1-e_3,e_2-e_3>$ and $f_{E_1}(x,y)=\frac{x}{y}$, where  $S=\mathbb{C}\setminus\{0,1\}$;\\
For the triple $L_3\cap L_4\cap L_5$, we have\\
$E_2^0=<e_4-e_5,e_4-e_3>$ and $f_{E_2}(x,y)=\frac{x-1}{y-1}$, where  $S=\mathbb{C}\setminus\{0,1\}$;\\
For the triple $L_2\cap L_4\cap L_5$, we have\\
$E_3^0=<e_4-e_2,e_6-e_2>$ and $f_{E_3}(x,y)=\frac{x-1}{y}$, where  $S=\mathbb{C}\setminus\{0,1\}$;\\
For the triple $L_1\cap L_5\cap L_7$, we have\\
$E_4^0=<e_1-e_7,e_5-e_1>$ and $f_{E_4}(x,y)=\frac{x}{y-1}$, where  $S=\mathbb{C}\setminus\{0,1\}$;\\
For the triple $L_1\cap L_4\cap L_8$, we have\\
$E_5^0=<e_1,e_4>$ and $f_{E_6}(x,y)=x$, where  $S=\mathbb{C}\setminus\{0,1\}$;\\
For the triple $L_2\cap L_5\cap L_8$, we have\\
$E_6^0=<e_5,e_2>$ and $f_{E_7}(x,y)=y$, where  $S=\mathbb{C}\setminus\{0,1\}$;\\
For the quadruple $L_3\cap L_6\cap L_7\cap L_8$, we have\\
$E_7^0=<e_3,e_6,e_7>$ and $f_{E_5}(x,y)=x-y$, where
$S=\mathbb{C}\setminus\{0,\pm1\}$.

\medskip

{2. \bf The non-local components} There are 5 non-local components,
corresponding to braid
subarrangements:\\
For $X=(L_1L_6|L_3L_4|L_2L_8)$, we have\\
$E_8^0=<e_1-e_3-e_4+e_6,e_2-e_3-e_4>$ and $f_{E_8}(x,y)=\frac{x(x-y-1)}{(x-y)(x-1)}$, where  $S=\mathbb{C}\setminus\{0,1\}$;\\
For $Y=(L_4L_8|L_2L_3|L_5L_6)$, we have\\
$E_9^0=<-e_2-e_3+e_5+e_6,e_2+e_3-e_4>$ and $f_{E_9}(x,y)=\frac{x-1}{y(x-y)}$, $S=\mathbb{C}\setminus\{0,1\}$\\
For $Z=(L_1L_5|L_2L_4|L_3L_8)$, we have\\
$E_{10}^0=<e_1-e_2-e_4+e_5,e_2-e_3+e_4>$ and $f_{E_{10}}(x,y)=\frac{x(y-1)}{y(x-1)}$, where  $S=\mathbb{C}\setminus\{0,1\}$;\\
For $W=(L_1L_3|L_4L_7|L_5L_8)$, we have\\
$E_{11}^0=<e_1+e_3-e_5,e_5-e_7-e_4>$ and
$f_{E_{11}}(x,y)=\frac{x(x-y)}{(x-y+1)(x-1)}$,
where  $S=\mathbb{C}\setminus\{0,1\}$;\\
For $V=(L_1L_8|L_2L_7|L_3L_5)$, we have\\
$E_{12}^0=<e_1-e_2-e_7,e_3+e_5-e_2-e_7>$ and
$f_{E_{12}}(x,y)=\frac{x}{y(x-y+1)}$,
where  $S=\mathbb{C}\setminus\{0,1\}$.\\

\medskip

One way to obtain these 12 irreducible components $E_j$ is to
compute the cup-product
$$H^1(M,\C)\times H^1(M,\C) \to H^2(M,\C)$$
and then to use the computer program SINGULAR to list the
irreducible components of the determinantal variety corresponding to
$R_1(d\A)$, see for details \cite{N}.

For each $f_E$ in the list above we can use the method described in
 \cite{D1} and we get that there is no translated
 component in $\V_1(d\A)$ associated to such an $f_E$.

 It was discovered by A. Suciu (again by using some computer computations) that $\V_1(d\A)$
 has one  1-dimensional translated component  $W$ associated to the
 mapping $f:M\rightarrow \C^*$ defined as, in affine coordinates,$$f(x,y)=\frac{x(y-1)(x-y-1)^2}{(x-1)y(x-y+1)^2}$$
and with $\rho_W=(1,-1,-1,-1,1,1,1)\in(\C^*)^7$, see  \cite{S1},
\cite{S2}. In other words, $W=\rho_W \otimes
\{(t,t^{-1},1,t^{-1},t,t^2,t^{-2},1)\ |\
 t\in \C^* \}$.

 Let $V_i$ be the component of $\V_1(d\A)$ corresponding to each
 $E_i$ for $i=1,...,12$, i.e. $V_i=\exp(E_i)$. Then it is known that
 $$W\cap V_{8}\cap V_9\cap V_{10}=\rho_W$$
 and
 $$W\cap V_{10}\cap V_{11}\cap V_{12}=\rho'_W,$$
 where $\rho'_W=(-1,1,-1,1,-1,1,1)\in(\C^*)^7$, see \cite{S1}, \cite{S2}.
 Since these results were obtained as a result of computer computations, it is useful to provide a
 direct proof based on Corollary \ref{cor1}.

 Let $A=E_8^0$ and $B=E_9^0$ be the primitive lattices in $H^1(M,\Z)$
 introduced above and apply to them the construction explained in Section 2
 with $L=H^1(M,\Z)$. Here $n=7$, $a=b=2$. The basis $e'_1$,... $e'_5$ can be chosen
 as given by the following equivalence classes
 $$ e'_1=[e_1-e_3-e_4+e_6],~~ e'_2=[e_3], ~~e'_3=[e_2],~~ e'_4=[e_5] ,~~e'_5=[e_7].$$
 Then $m=1$ and $d_2=2$. Let $f_1=e_2+e_3-e_5-e_6$ and $f_2=e_4-e_2-e_3$.
 Then $\B=\{e_1-e_3-e_4+e_6, e_3, e_2,e_5,e_7,f_1,f_2\}$ is a $\Z$-basis of $L$
 (the coefficient matrix is unimodular) and we can take $g_1=e_1-e_3-e_4+e_6$ and $g_2=2e_3+f_2$.
  Therefore $$Tors(H^1(M,\Z)/(E_8^0+E_9^0))=\Z_2.$$
  Now, by the morphism
$\theta:\Z_2\rightarrow V_8\cap V_9$ used in
 section \ref{sec:proof},
  $$\hat{1}\mapsto\exp_L(\frac{1}{2}(g_2))=\exp_L(\frac{1}{2}(-e_2+e_3+e_4))=(1,-1,-1,-1,1,1,1) =\rho_W.
 $$
By Corollary \ref{cor1}, it follows that
$$V_8\cap V_9=\theta(\Z_2)=\{1,\rho_W\}.$$
In exactly the same way one can show that $V_8\cap V_{10}=V_9 \cap
V_{10}=\{1,\rho_W\}.$ Since clearly $\rho_W\in W$, it follows that
the four components $V_8$, $V_9$, $ V_{10}$ and $W$ meet exactly in
one point.

Similarly one shows that $W\cap V_{10}\cap V_{11}\cap
V_{12}=\rho'_W$ and that all the other intersections of two
irreducible components are trivial.

\end{document}